# Finite Projective Planes

Dhananjay P. Mehendale
Sir Parashurambhau College, Tilak Road, Pune 411030,
India

## Abstract

We propose graph theoretic equivalents for existence of a finite projective plane (*fpp*). We then develop a new approach and see that the problem of existence of a finite projective plane of order *n* is linked up with a subset containing $n(n-1)$ sharply 2-transitive permutations. If *n* is prime power then it is well known that there exists a finite field and existence of this field implies existence of (*n*-1) MOLS which further implies existence of *fpp*. We show that by assuming the existence of (*n*-1) MOLS the existence of a group of order $n(n-1)$ made up of sharply 2-transitive permutations can be implied through transforming the given (*n*-1) MOLS to suitable form. From a known result [1], [2] it then follows that when such group exists the order *n* has to be a prime power. Finally, we see relation between MOLS and determinantal monomials and between MOLS and a cyclic group that permutes the rows of MOLS. Finally, we conclude the paper with some important remarks.

1. **Introduction:** When *n* is prime power one can suitably define binary operations + and * modulo *n* and construct finite field of *n* elements, $\{0,1,2,\cdots,n-1\}$, called Galois field, GF(*n*), modulo *n*. Using this Galois field it is possible to construct (*n*-1) Mutually Orthogonal Latin Squares (MOLS) using standard (formula for (*i*, *j*)-th element in *k*-th Latin Square and procedure for their) construction. Thus, there exists finite field of order *n* when $n = p^k$, where *p* is prime and *k* is some positive integer greater than one, and using this field one can construct (*n*-1) MOLS made up of elements of the finite field under consideration as the matrix elements. It is a well known result due to R. C. Bose [5] that there exist (*n*-1) MOLS if and only if a finite projective plane of order *n* exists. Therefore, when *n* is prime power finite projective plane of order *n* exists. The main question associated with finite projective planes (*fpps*) is whether there exist finite projective planes of nonprime power orders. All the so far known *fpps* are of prime power order and it has been conjectured that the order of an *fpp* is always either a prime or a prime power. A straight forward way to settle this conjecture is to show nonexistence of *fpp* for all nonprime power orders. The following celebrated result due to Bruck and Ryser [3] establishes the nonexistence of *fpp* for certain types of infinitely many (but not all) nonprime power orders.

**Theorem 1.1 ([3]):** If $n \equiv 1$ or $2 \mod (4)$ and if the square free part of $n$ contains at least one prime factor of the form $4k+3$ then there does not exist any projective plane of order $n$.

$\square$

This theorem establishes the nonexistence of projective plane for infinitely many values as orders, like, $n = 6, 14, 21, \cdots$ however, it leaves behind equally many values as orders, like, $n = 10, 12, 15, \cdots$ for which the theorem is inconclusive and cannot decide whether or not the *fpp* for these orders exists.

The case $n = 10$ was successfully resolved making heavy use of computers by C. W. H. Lam [4].

2. **Graph Theoretic Equivalent for *fpp*:** We know that a finite projective plane is by definition a collection of $n^2 + n + 1$ points and $n^2 + n + 1$ lines such that the following four axioms hold:
   (i) Every line contains $n + 1$ points.
   (ii) Every point lies on $n + 1$ lines.
   (iii) Any two distinct lines intersect in exactly one point, and
   (iv) Any two distinct points lie on exactly one line.
   In the language of design theory a finite projective plane is nothing but a $2 - (n^2 + n + 1, n + 1, 1)-$design, or a $S(n^2 + n + 1, n + 1, 1)$ Steiner system.

**Definition 2.1:** A graph is called **monochromatic** if all its edges have same color.

We now state a graph theoretic equivalent of the problem of existence of *fpp*:

**Definition 2.2:** A set of vertices is called **totally disconnected** if no two vertices of this set are joined by an edge.

**Definition 2.3:** A **totally disconnected graph** is a graph whose vertex set is a totally disconnected set.

**Theorem 2.1:** Suppose we are given ($n^2 + n + 1$) monochromatic complete graphs on ($n+1$) points such that each one is colored in some unique color different from the one used for any other complete graph then we can **tightly pack** these complete graphs on ($n+1$) points in a totally disconnected graph on ($n^2 + n + 1$) points if and only if the projective plane of order $n$ exists.

**Proof:** Suppose the tight packing exists. We write down the vertex labels of these complete graphs on ($n+1$) points in increasing order in rows, forming in total $n^2 + n + 1$ rows. Each such row will represent vertices of a complete graph on ($n+1$) vertices. It is straightforward to see that these rows will constitute the $2-(n^2+n+1, n+1, 1)$−design. Conversely, if projective plane of order $n$ exists then the $2-(n^2+n+1, n+1, 1)$−design exists and using the numbers in its rows as labels of vertices for the complete graphs on ($n+1$) points we can attain the desired tight packing.

**Illustrative Examples:**   □

The first figure given below represents packing of 7 monochromatic triangles in $K_7$.

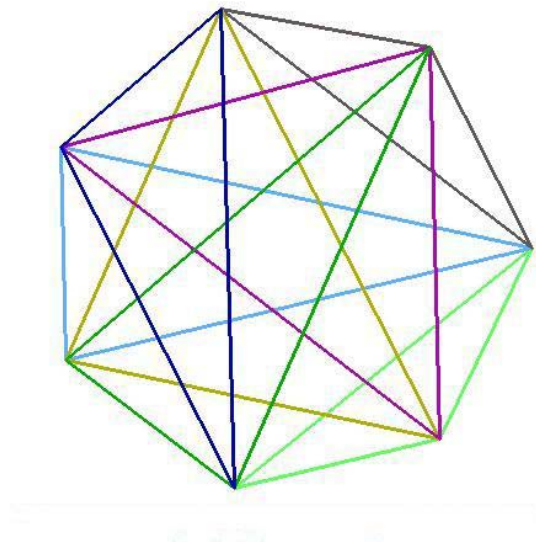

The second figure given below represents packing of 13 monochromatic $K_4$ in $K_{13}$.

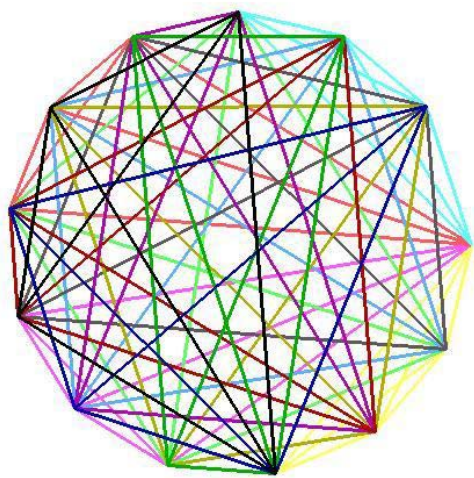

Let the points (vertices) of first graph (G) and second graph (H) above be labeled as V(G) = {1, 2, ….,7} and V(H) = {1, 2, …13}. The rows of vertices forming lines (monochromatic complete graphs) for these cases are

```
1 2 3
1 4 5
1 6 7
2 4 6
2 5 7
3 4 7
3 5 6
```

and

```
1  2  3  4
1  5  6  7
1  8  9  10
1  11 12 13
2  5  8  11
2  6  9  12
```

```
2   7   10   13
3   5    9   13
3   6   10   11
3   7    8   12
4   5   10   12
4   6    8   13
4   7    9   11
```

From the above given 7 subsets each containing 3 vertices (respectively 13 subsets each containing 4 vertices) made out of set of 7 vertices (respectively 13 vertices) the most important property they should posses to form a projective plane of order 2 (respectively order 3) is as follows: Any two of these 3-subsets (respectively 4-subsets) of vertices should have at most one vertex in common. As previous, a subset of vertices is called **totally** disconnected if no two vertices of this subset are joined by an edge. The tight packing implies the existence of 7 totally disconnected 3-subsets of vertices formed out of 7 vertices (respectively, 13 totally disconnected 4-subsets of vertices formed out of 13 vertices). How many equivalent representations are possible like the one given above in these two cases? To know and determine such representations we may proceed with so called **"edge labeling algorithm"** as follows:

(i) Start with totally disconnected graph of seven (respectively thirteen) vertices.
(ii) Take any totally disconnected set of 3 vertices (respectively 4 vertices) out of this totally disconnected set of 7 vertices (respectively 13 vertices) and join all its edges and assign color "1" to all these edges.
(iii) Take another totally disconnected set of 3 vertices (respectively 4 vertices) and join all its edges and assign next color, say "2" to all these edges.
(iv) Continue edge labeling in this fashion till finally every pair of vertices gets joined by an edge of some color and thus every edge will get some unique color.
(v) Now assign labeling to vertices in any way and collect the 3-sets (respectively 4-sets) of thus labeled vertices whose edges (joining them to each other) are colored by same color.

It is now clear to see that when a projective plane of some order $n$ doesn't exist then using chosen totally disconnected set of $(n^2 + n + 1)$ vertices one can't form ($n^2 + n + 1$) totally disconnected (*n*+1)-subsets of vertices and therefore edge labeling algorithm mentioned above cannot be carried out which requires existence of such subsets. From this graph theoretic equivalent of *fpp* problem it is clear that the existence of projective plane of order *n* is equivalent to existence of a complete graph on ($n^2 + n + 1$) vertices, containing ($n^2 + n + 1$) monochromatic complete subgraphs each containing ($n + 1$) vertices with edges colored by some chosen color, and which are different for each such different complete subgraph. Due to this graph theoretic equivalent, we have following associations: Here, the ($n^2 + n + 1$) vertices are "points" while complete subgraphs each containing ($n + 1$) vertices are "lines" for our **Model of finite projective plane**. It is easy to check that all the defining properties of projective plane are satisfied by this model, i.e. a finite projective plane is collection of ($n^2 + n + 1$) vertices and ($n^2 + n + 1$) monochromatic complete subgraphs each containing ($n + 1$) vertices, and whose edges are colored with some chosen color for each ($n + 1$)-subgraph, such that the following four axioms hold:

(i) Every monochromatic complete subgraph contains ($n + 1$) vertices.
(ii) Every vertex lies on ($n + 1$) monochromatic complete subgraphs.
(iii) Any two distinct monochromatic complete subgraphs intersect in exactly one vertex, and
(iv) Any two distinct vertices lie on exactly one monochromatic complete subgraph.

Considering the above axioms together suggests that the following **steps** are advised while trying to find a **tight packing** mentioned above for some order *n*.
 **Step 1)** Label the vertices with numbers {1, 2, 3, ……, $n^2 + n + 1$}
 **Step 2)** Form complete subgraphs on ($n + 1$) vertices such that exactly ($n + 1$) complete subgraphs will pass through vertex with label 1 and only 1 complete subgraph will pass through vertices with label {2, 3, 4, ….., $n^2 + n + 1$}.
 **Step 3)** Form complete subgraphs on ($n + 1$) vertices such that exactly ($n + 1$) complete subgraphs will pass through vertex with label {1, 2} and only 1 complete subgraph will pass through vertices with label {3, 4, …, $n + 1$} and exactly 2 complete subgraphs will pass through vertices {*n*+2, *n*+3, ….., $n^2 + n + 1$}

**Step 4)** Form complete subgraphs on $(n+1)$ vertices such that exactly $(n+1)$ complete subgraphs will pass through vertex with label {1, 2, 3} and only 1 complete subgraph will pass through vertices with label {4, 5, ..., $n+1$} and exactly 3 complete subgraphs will pass through vertices {$n+2, n+3, ....., n^2 + n + 1$}

 **Step 5)** Continue these steps till finally we (can) form complete subgraphs on $(n+1)$ vertices such that exactly $(n+1)$ complete subgraphs will pass through vertex with label {1, 2, 3, ....., $n+1$} and exactly $(n+1)$ complete subgraphs will pass through vertices {$n+2, n+3, ....., n^2 + n + 1$}. Thus, finally exactly $(n+1)$ complete subgraphs will pass through ALL vertices!

Another graph theoretic equivalent for the problem of existence of *fpp* can be as given below:

A permutation $P$ on $n$ symbols expressed in terms of disjoint cycles can be viewed as collection of **directed disjoint cycles** in a directed graph on $n$ points, e.g. Let $P$ = (14)(25)(36) then $P$ can be viewed as collection of three disjoint directed cycles in a digraph on six points and made up of cycles $1 \to 4 \to 1$, $2 \to 5 \to 2$, $3 \to 6 \to 3$.

**Definition 2.3:** A **complete cycle cover** for a symmetric digraph on $n$ points is a collection of permutations of degree $n$ whose disjoint cycle representation when viewed as collection of disjoint directed cycles in a digraph on $n$ points together exactly cover the entire symmetric digraph, i.e. these disjoint cycles together tightly pack the symmetric digraph on $n$ points such that every directed edge belongs to some disjoint directed cycle.

**Theorem 2.2:** A complete cycle cover formed from a group of $n(n-1)$ permutations on $n(n-1)$ symbols for a symmetric digraph on $n(n-1)$ points exists if and only if the projective plane of order $n$ exists.

**Proof:** Consider the group of $n(n-1)$ permutations on $n(n-1)$ symbols forming the complete cycle cover for the symmetric digraph on $n(n-1)$ points. We can think of this group as **pair group** of the group of degree $n$ and of order $n(n-1)$. This group of degree $n$ and of order $n(n-1)$ can be obtained from mapping: $1 \to (1,2)$, $2 \to (1,3)$, ...., $n(n-1) \to (n,(n-1))$. One can check that this recovered original group is actually a 2-transitive group of degree $n$ and order $n(n-1)$ containing a regular normal subgroup of order $n$. This will imply the existence of a projective plane as per the discussion in section 5. Conversely, when the projective plane of order $n$ exists we will have the existence of a 2-transitive subgroup of $S_n$ of order $n(n-1)$ containing a regular normal subgroup of order $n$. Doing the reverse mapping we can construct the pair group of this group of both the degree and order $n(n-1)$ and show the existence of the desired complete cycle cover.

☐

**Illustrative Example:**

| Cycle Representation for group $S_3$ | Cycle Representation for its Pair Group |
|---|---|
| (1)(2)(3) | (1)(2)(3)(4)(5)(6) |
| (123) | (135)(246) |
| (132) | (153)(264) |
| (1)(23) | (12)(36)(45) |
| (2)(13) | (16)(25)(34) |
| (3)(12) | (14)(23)(56) |

The complete cycle cover for the symmetric digraph on six points using the cycle representation for the pair group given above is made up of the following directed disjoint cycles:

(1) $1 \to 1, 2 \to 2, 3 \to 3, 4 \to 4, 5 \to 5, 6 \to 6$.
(2) $1 \to 3 \to 5 \to 1, 2 \to 4 \to 6 \to 2$.
(3) $1 \to 5 \to 3 \to 1, 2 \to 6 \to 4 \to 2$.
(4) $1 \to 2 \to 1, 3 \to 6 \to 3, 4 \to 5 \to 4$.
(5) $1 \to 6 \to 1, 2 \to 5 \to 2, 3 \to 4 \to 3$.
(6) $1 \to 4 \to 1, 2 \to 3 \to 2, 5 \to 6 \to 5$.

3. **A New Approach:** In our new approach we suggest to build the projective plane of order $n$ as follows:

(1) We first build a set of $(2n+1)$ lines of the projective plane. This set contains a subset of $(n+1)$ lines and we call these lines in this subset the **horizontal lines** of the projective plane and a subset of $n$ lines and we call these lines in this subset the **vertical lines** of the projective plane.
(2) We then proceed, if and when possible, with the construction of the subset of remaining $(n^2 - n)$ lines and we call these lines in this subset the **determinantal lines** of the projective plane. When we are successful with this construction then we have successfully completed the construction of the desired finite projective plane of order $n$.

The first step mentioned above is always possible to achieve for any value of $n$. <u>But, what are the values of $n$ for which the second step is also possible?</u> We want to investigate whether the second step is possible only when $n$ is prime or prime power.
Let us begin with the first step of building $(2n+1)$ lines. The first step of building $(2n+1)$ lines begins with the construction of so called $(n+1)$ **horizontal lines** of the projective plane and ends with the construction of remaining $n$ lines the so called $n$ **vertical lines** of the projective plane.

The **horizontal lines** of the projective plane are:

$$
\begin{array}{lllll}
1 & 2 & 3 & \cdots & (n+1) \\
1 & (n+2) & (n+3) & \cdots & (2n+1) \\
1 & (2n+2) & (2n+3) & \cdots & (3n+1) \\
\vdots & \vdots & \vdots & \ldots & \vdots \\
1 & (n^2+2) & (n^2+3) & \cdots & (n^2+n+1)
\end{array}
\qquad (2.1)
$$

The **vertical lines** of the projective plane are:

$$\begin{array}{ccccc} 2 & (n+2) & (2n+2) & \cdots & (n^2+2) \\ 2 & (n+3) & (2n+3) & \cdots & (n^2+3) \\ 2 & (n+4) & (2n+4) & \cdots & (n^2+4) \\ \vdots & \vdots & \vdots & \cdots & \vdots \\ 2 & (2n+1) & (3n+1) & \cdots & (n^2+n+1) \end{array} \qquad (2.2)$$

The second step consists of building the remaining $(n^2 - n)$ so called **determinantal lines** of the projective plane (whenever possible). We know that whenever $n$ is prime or prime power the projective plane for that order exists and so in this case it must be always possible to build the desired $(n^2 - n)$ determinantal lines of the projective plane along with $(n+1)$ horizontal lines and $n$ vertical lines. Our aim is to see that whenever we can also build the desired $(n^2 - n)$ determinantal lines along with $(n+1)$ horizontal lines and $n$ vertical lines then why $n$ has to be **prime power** in such case. The determinantal lines are made up of $(n-1)$ subsets of lines taken together, such that each of these subsets contains $n$ lines. The leading entry of the lines in these subsets, respectively the first, second, ...., $(n-1)$-th subsets of these lines, are $3, 4, \cdots, (n+1)$. The other entries in these lines are taken by selecting appropriate entries from the sub-table formed by deleting first row and first column of the table of horizontal lines (2.1), such that exactly one entry from each row and column (like a determinantal monomial) of the sub-table gets incorporated in each of the so called determinantal line.

In order to illustrate this new approach we proceed with

**Few Examples:** As per the procedure described above we begin with the construction of some projective planes of small orders:

(1) The well-known **Feno plane,** the *fpp* of order $n = 2$:

Horizontal lines:
1 2 3
1 4 5
1 6 7

Vertical lines:
2 4 6
2 5 7

Determinantal lines:
3 4 7
3 5 6

Now, we can think about the entries (4 7) in the above given first determinantal line as arrived at from permutation $\begin{pmatrix} 1 & 2 \\ 1 & 2 \end{pmatrix} = (1)(2)$, while the entries (5 6) in the second determinantal line as arrived at from permutation $\begin{pmatrix} 1 & 2 \\ 2 & 1 \end{pmatrix} = (12)$. Thus, we can think of the determinantal lines as generated by the group $S_2$, a sharply 2-transitive group of permutations on 2 symbols containing the only subgroup of order one formed by identity element.

Now, we note down the parts of determinantal lines (arrived as determinantal monomials of the matrix formed by the sub-table resulting from deleting the first row and first column of the table of horizontal lines) given above and visualize them as permutations. We write in front of these parts the corresponding cycle decomposition and transformed row representation of the permutation under consideration:

| Parts | Cycle decomposition | Transformed Rows |
|---|---|---|
| 4  7 | (1)(2) | 1  2 |
| 5  6 | (12) | 2  1 |

We can split the $(n^2 - n)$ transformed rows representing permutations into $(n-1)$ subsets such that each of these subsets contains $n$ transformed rows and build square matrices of size $n$, say $L^c + 1$, where $L$ represents the Latin Square. Since $(n-1) = 1$, therefore in the present case we have

$$L^c + 1 = \begin{matrix} 1 & 2 \\ 2 & 1 \end{matrix}$$

therefore, we get only one Latin Square:

$$L = \begin{matrix} 0 & 1 \\ 1 & 0 \end{matrix}$$

Further, it can be easily checked that by using cycle decomposition corresponding to each transformed row among the rows forming an $L^c + 1$ and superimposing the directed graphs corresponding to each transformed row representing a permutation we can build a cycle cover that forms a complete symmetric diagraph. Thus each $L^c + 1$ gives rise to a (different) cycle cover which when superimposed leads to formation of a complete symmetric diagraph. The following is the cycle cover leading to formation of complete symmetric diagraph on two vertices:

$$1 \to 1, 2 \to 2$$
$$1 \to 2 \to 1$$

(2) The *fpp* of order $n = 3$:

Horizontal lines:
```
1   2   3   4
1   5   6   7
1   8   9  10
1  11  12  13
```

Vertical lines:
2   5   8   11
2   6   9   12
2   7   10  13

Determinantal lines:
3   5   9   13
3   6   10  11
3   7   8   12

4   5   10  12
4   6   8   13
4   7   9   11

A permutation of degree $n$ can be expressed in the standard notation as $\begin{pmatrix} 1 & 2 & \cdots & n-1 & n \\ i_1 & i_2 & \cdots & i_{(n-1)} & i_n \end{pmatrix}$. Let us call the lower row, $i_1 \; i_2 \; \cdots \; i_{(n-1)} \; i_n$, in the above representation for a permutation the **transformed row** of numbers due to the action of the permutation which takes $1 \to i_1, 2 \to i_2, \cdots$ etc..

| Parts | Cycle decomposition | Transformed Rows |
|---|---|---|
| 5   9   13 | (1)(2)(3) | 1 2 3 |
| 6  10   11 | (123)     | 2 3 1 |
| 7   8   12 | (132)     | 3 1 2 |
|            |           |       |
| 5  10   12 | (1)(23)   | 1 3 2 |
| 6   8   13 | (3)(12)   | 2 1 3 |
| 7   9   11 | (2)(13)   | 3 2 1 |

Clearly, we can think of the determinantal lines as generated by the sharply 2-transitive group $S_3$ the group of permutations on 3 symbols, containing a regularly normal subgroup of order 3 formed by two fixed point free permutations $\{(123), (132)\}$ and the identity $\{(1)(2)(3)\}$, i.e. corresponding to first three permutations given above.

Thus, for the formation of finite projective plane existence of $(n^2 - n) = n(n-1)$ determinantal lines is important and for existence of determinantal lines the existence of $(n^2 - n) = n(n-1)$ sharply 2-transitive permutations is important. So, we now pick the second and third columns given above and note below six sharply 2-transitive permutations forming group $S_3$ in a two-columned table. The first column gives the disjoint cycle decomposition of these permutations and second column gives the corresponding transformed rows for these permutations. All these permutations together form the sharply 2-transitive group $S_3$:

| Permutation(disjoint cycle form) | Transformed rows |
|---|---|
| (1)(2)(3) | 1 2 3 |
| (123) | 2 3 1 |
| (132) | 3 1 2 |
| (1)(23) | 1 3 2 |
| (2)(13) | 3 2 1 |
| (3)(12) | 2 1 3 |

**Remark 3.1:** Since, by the very nature of their construction procedure, we can always build the horizontal and the vertical lines therefore the thing of real importance is the **existence of** $(n^2 - n) = n(n-1)$ **determinantal lines** for the completion of the construction of the desired projective plane. For the existence of $(n^2 - n) = n(n-1)$ determinantal lines existence of $(n^2 - n) = n(n-1)$ sharply 2-transitive permutations is important. In the case of the above examples the groups $S_2$, $S_3$ themselves proved useful for the construction of desired *fpp* but we will see below that for orders $n \geq 4$ we need to see whether there exists a sharply 2-transitive subgroup of $S_n$ of order $n(n-1)$ containing a regularly normal subgroup of order $n$ formed by **fixed-point-free permutations** and the identity element. Thus, in order to show the existence of *fpp* of order $n$ we should try to see whether we can find a sharply 2-transitive subgroup of $S_n$, of order $n(n-1)$ containing a regularly normal subgroup of order $n$.

In the present example $S_3$ is sharply 2-transitive group having regularly normal subgroup of order 3, containing elements {(1)(2)(3), (123), (132)}, i.e. two fixed point free permutations and the identity.

It is easy to check that in the present example one can construct two $L^c + 1$ using the transformed rows given in the second column of the above table which leads to formation of two MOLS as follows:

$$L_1^c + 1 = \begin{matrix} 1 & 2 & 3 \\ 2 & 3 & 1 \\ 3 & 1 & 2 \end{matrix}, \quad L_2^c + 1 = \begin{matrix} 1 & 3 & 2 \\ 2 & 1 & 3 \\ 3 & 2 & 1 \end{matrix}$$

$$L_1 = \begin{matrix} 0 & 1 & 2 \\ 1 & 2 & 0 \\ 2 & 0 & 1 \end{matrix}, \quad L_2 = \begin{matrix} 0 & 1 & 2 \\ 2 & 0 & 1 \\ 1 & 2 & 0 \end{matrix}$$

The cycle decomposition corresponding to transformed rows for $L_1^C + 1$ is (1)(2)(3)   and for $L_2^C + 1$ is (1)(23)
 (123)                 (12)(3)
 (132)                 (13)(2)

So clearly, the cycle decompositions corresponding to all transformed rows taken together give rise to $S_3$, a sharply 2-transitive group having regularly normal subgroup of order 3 containing two 3-cycles and identity. Further, cycle decomposition leads to formation of complete symmetric diagraph. Two Cycle covers leading to formation of complete symmetric diagraph on three vertices in two different ways occurs as below. First cycle cover is formed using transformed rows in $L_1^C + 1$ and second cycle cover is formed using transformed rows in $L_2^C + 1$:

$$1 \to 1, 2 \to 2, 3 \to 3$$
$$1 \to 2 \to 3 \to 1, 1 \to 3 \to 2 \to 1$$

and

$$1 \to 1, 2 \to 3 \to 2$$
$$3 \to 3, 1 \to 2 \to 1$$
$$2 \to 2, 1 \to 3 \to 1$$

(3) The *fpp* of order $n = 4$:

The suitable group of order $n(n-1) = 12$ for this case is $A_4$, the alternating group, this group $A_4$ is sharply 2-transitive subgroup of $S_4$ made up of even permutations and having normal subgroup of order 4, containing elements {(1)(2)(3)(4), (12)(34), (13)(24), (14)(23)}, i.e. three fixed point free permutations and the identity.

A permutation of degree $n$ can be expressed in the standard notation as $\begin{pmatrix} 1 & 2 & \cdots & n-1 & n \\ i_1 & i_2 & \cdots & i_{(n-1)} & i_n \end{pmatrix}$. Let us call the lower row, $i_1 \; i_2 \; \cdots \; i_{(n-1)} \; i_n$, in the above representation for a permutation the **transformed row** of numbers due to the action of the permutation which takes $1 \to i_1, 2 \to i_2, \cdots$ etc..

We now note below in a two-columned table, the first column gives the disjoint cycle decomposition of the respective permutation and second column gives the corresponding transformed row for the same permutation. All these permutations together form the group $A_4$, a sharply 2-transitive subgroup of $S_4$:

| Permutation(disjoint cycle form) | Transformed rows |
|---|---|
| (1)(2)(3)(4) | 1 2 3 4 |
| (1)(234) | 1 3 4 2 |
| (1)(243) | 1 4 2 3 |
| (12)(34) | 2 1 4 3 |
| (123)(4) | 2 3 1 4 |
| (124)(3) | 2 4 3 1 |
| (132)(4) | 3 1 2 4 |
| (134)(2) | 3 2 4 1 |

| (13)(24) | 3 4 1 2 |
| (142)(3) | 4 1 3 2 |
| (143)(2) | 4 2 1 3 |
| (14)(23) | 4 3 2 1 |

It is easy to check that one can construct $L_1^C +1$, $L_2^C +1$, $L_3^C +1$ using the transformed rows given in the second column of the above table as follows:

| **1 2 3 4** | 1 3 4 2 | 1 4 2 3 |
| **2 1 4 3** | 2 4 3 1 | 2 3 1 4 |
| **3 4 1 2** | 3 1 2 4 | 3 2 4 1 |
| **4 3 2 1** | 4 2 1 3 | 4 1 3 2 |

The permutations corresponding to transformed rows of $L_1^C +1$ form the normal subgroup. In this case, the transformed rows corresponding to the normal subgroup of order 4 mentioned above. Note that there are exactly $n(n-1) = 12$ transformed rows from which we can construct ($n-1=3$) MOLS, each containing $n=4$ rows (and columns), which implies the existence *fpp*. It is clear to see that $\{L_1^C +1, L_2^C +1, L_3^C +1\}$ can be treated as MOLS and thus, rows MOLS can be thought of as made up of (transformed rows corresponding to) sharply 2-transitive permutations. The transformed rows of $L_1^C +1$ forms normal subgroup of group of sharply 2-transitive permutations formed by transformed rows of all $\{L_1^C +1, L_2^C +1, L_3^C +1\}$ taken together.

Further, three Cycle covers leading to formation of complete symmetric diagraph on four vertices in three different ways for $L_1^C +1$ and $L_2^C +1$ and $L_3^C +1$ are respectively:

$$1 \to 1, 2 \to 2, 3 \to 3, 4 \to 4$$
$$1 \to 2 \to 1, 3 \to 4 \to 3$$
$$1 \to 3 \to 1, 2 \to 4 \to 2$$
$$1 \to 4 \to 1, 2 \to 3 \to 2$$

and

$$1 \to 1, 2 \to 3 \to 4 \to 2$$
$$1 \to 2 \to 4 \to 1, 3 \to 3$$
$$1 \to 3 \to 2 \to 1, 4 \to 4$$
$$1 \to 4 \to 3 \to 1, 2 \to 2$$

and

$$1 \to 1, 2 \to 4 \to 3 \to 2$$
$$1 \to 2 \to 3 \to 1, 4 \to 4$$
$$1 \to 3 \to 4 \to 1, 2 \to 2$$
$$1 \to 4 \to 2 \to 1, 3 \to 3$$

(4) The *fpp* of order $n = 5$:

The suitable group of order $n(n-1) = 20$ for this case is the following (*k*-metacyclic) sharply 2-transitive subgroup $G$ of group $S_5$, satisfying the generator relations $s^5 = t^4 = e$, the identity, and $t^{-1}st = s^2$, for all $s, t \in G$, and containing a normal subgroup of order 5, formed by four elements (permutations) without containing a fixed point and the identity, i.e. containing four 5-cycles and the identity element.

We now give below a two-columned table for this (*k*-metacyclic) group, $G$. The disjoint cycle representation for permutations which are elements of this group is given in the first column and the transformed rows corresponding to these permutations are given in the same rows in the second column of this two-columned table given below:

| Permutation(disjoint cycle form) | Transformed rows |
| --- | --- |
| (1)(2)(3)(4)(5) | 1 2 3 4 5 |
| (1)(2354) | 1 3 5 2 4 |
| (1)(2453) | 1 4 2 5 3 |
| (1)(25)(34) | 1 5 4 3 2 |
| (12345) | 2 3 4 5 1 |
| (1243)(5) | 2 4 1 3 5 |
| (1254)(3) | 2 5 3 1 4 |
| (12)(35)(4) | 2 1 5 4 3 |
| (13524) | 3 4 5 1 2 |
| (1325)(4) | 3 5 2 4 1 |
| (1342)(5) | 3 1 4 2 5 |
| (13)(2)(45) | 3 2 1 5 4 |
| (14253) | 4 5 1 2 3 |
| (1452)(3) | 4 1 3 5 2 |
| (1435)(2) | 4 2 5 3 1 |
| (14)(23)(5) | 4 3 2 1 5 |
| (15432) | 5 1 2 3 4 |
| (1534)(2) | 5 2 4 1 3 |
| (1523)(4) | 5 3 1 4 2 |
| (15)(24)(3) | 5 4 3 2 1 |

Now, it is easy to see as follows that we can construct four $L^c + 1$ (which lead to formation of four MOLS) using the transformed rows given in the second column of the above table. The first

$L^c + 1$ contains transformed rows corresponding to the elements of the normal subgroup of order 5 formed by identity and four fixed point free permutations:

```
1 2 3 4 5     1 3 5 2 4     1 4 2 5 3     1 5 4 3 2
2 3 4 5 1     2 4 1 3 5     2 5 3 1 4     2 1 5 4 3
3 4 5 1 2     3 5 2 4 1     3 1 4 5 2     3 2 1 5 4
4 5 1 2 3     4 1 3 5 2     4 2 5 3 1     4 3 2 1 5
5 1 2 3 4     5 2 4 1 3     5 3 1 4 2     5 4 3 2 1
```

Further, Four Cycle covers such that each cycle cover is in fact different and is leading to formation of complete symmetric diagraph on five vertices in four different ways can be formed as is done in the previous cases.

**Specialty of these Groups:** As already mentioned, these groups of degree $n$ and order $n(n-1)$ are made up of **sharply 2-transitive** permutations.

Note that in the two columned tables given above we have in the second column the transformed rows written one below the other. Since transformed rows are written one below the other and each such row contains exactly $n$ entries so if we consider together the first entry in each row then we can say that they together form first column. Similarly, all the second entries together will form second column and so on. Thus, these $n(n-1)$ transformed rows placed one below the other such that each row is containing $n$ entries gives rise in effect to $n$ columns. Sharp 2-transitivity ensures that if we consider any two columns out of these $n$ columns and form pairs of numbers by taking first entry in these pairs from first chosen column and second entry in these pairs from second chosen column then there is exactly one time appearance of every pair among the pairs (1, 2), (1, 3), **….**, (1, n), (2, 1), (2, 3), **….**, (2, n), (3, 1), (3,2), (3,4), **….**, ((n-1), n), (n ,1), (n, 2), **….**, (n, (n-1)) in each such pair of columns $(i, j)$ such that $1 \leq i \leq n$ and $1 \leq j \leq n$. Exactly one time appearance of "**each pair of numbers $(i, j)$ in each column pair $(k, l)$**" tells us that there exists exactly one (i.e. sharp) permutation in the group taking pair $(k, l) \to (i, j)$, (i.e. 2-transitivity).Thus, Any two columns of the sub-table of transformed rows should **contain** every pair (1, 2), (1, 3), **……**,
(1, n), (2, 1), (2, 3), **….**, (2, n), (3, 1), (3,2), (3,4), **….**, ((n-1), n), (n ,1), (n, 2), **….**, (n, (n-1)), and only once. Also, in terms of the conditions that should be satisfied by cycle decomposition of the permutations given above in rhe table in order to avoid the repetition and/or omission of above mentioned pairs (1, 2), (1, 3), **….**,etc.:

There should be **no repetition** of
(1) A pair $(i)(j)$
(2) A 2-cycle $(i, j)$

(3) A string of same sequence of numbers like $(\cdots i_1, i_2, \cdots, i_l, \cdots)$ of length $l \geq 3$ in different cycles.

(4) There should be no two or more permutations whose cycle decomposition contains same two successive numbers, like, $\cdots i_1, i_2 \cdots$ and $\cdots j_1, j_2 \cdots$. This implies that there exist two permutations which maps $i_1 \to i_2$ and $j_1 \to j_2$ which will further implies that in the column pair $(i_1, j_1)$ the number pair $(i_2, j_2)$ appears twice or more times. This hampers sharp 2-transitivity.

As per our suggested approach we have seen that the existence of a projective plane requires existence of $n(n-1)$ determinantal lines, and this further gives rise to the requirement of existence of subset of $n(n-1)$ permutations which are sharply 2-transitive. In the cases so far seen, where the order $n$ was 2, 3, 4, 5 respectively, we got such $n(n-1)$ subsets of sharply 2-transitive permutations. In fact as seen above when $n$ was 2, 3, 4, and 5 respectively these subsets of sharply 2-transitive permutations were **sharply 2-transitive subgroups** $S_2$, $S_3$, $A_4$, and (the k-metacyclic group) $G$.

Let us fix up some notation: Let $X$ denotes a finite set and $|X|$ denotes degree (cardinality) of $X$. A permutation of $X$ is a bijection $\alpha : X \to X$. Two such permutations $\alpha, \beta$ can be composed to give another permutation $\alpha\beta : X \to X$ defined by the rule $\alpha\beta(x) = \alpha(\beta(x))$. Under this operation of composition the set of permutations of $X$ form a group, Sym($X$), the symmetric group of permutations on $X$. If we choose the set $X = \{1, 2, \cdots, n\}$ then we write $S_n$ for Sym($X$) and we have the cardinality of $S_n$, $|S_n| = n!$. If $G$ is subgroup of Sym($X$), then we shall say that the pair $(G, X)$ is a permutation group of degree $|X|$ and $G$ acts on $X$. Now let us note down the following important theorem given in ([2], page 127). This theorem tells us that when we can find a sharply 2-transitive subgroup of $S_n$, of order $n(n-1)$ containing a regularly normal subgroup of order $n$ then such $n$ must be a prime power.

**Theorem 3.1([2], page 127):** If $(G, X)$ is a sharply 2-transitive group then $|X|$ is prime power, $p^r$. The fixed point free elements of $G$, together with identity, form a regular normal subgroup of G which is also a Sylow p-subgroup of G.

The following equivalents of above theorem can be found in ([1], page 57) or ([2], page 19):

**Theorem 3.2([1], page 57):** Let $G$ be $m$-transitive on $X$ and $H$ a regularly normal subgroup of $G$.
(1) If $m = 2$ then $n$ is power of some prime $p$ and $H$ is an elementary abelian $p$-group.
(2) If $m = 3$ then $n$ is power of 2 or $n = 3$ and $G = \text{Sym}(X)$.
(3) If $m \geq 4$ then $m = 4 = n$ and $G = \text{Sym}(X)$.

□

**Theorem 3.3([2], page 19):** If $(G, X)$ is $k$-transitive and has a regular normal subgroup $N$, then
(1) $k = 2$ implies $N \approx (Z_p)^n, |X| = p^n$, $p$ prime.
(2) $k = 3$ implies $N \approx (Z_2)^n$ or $N \approx Z_3$, $|X| = 2^n$ or 3.
(3) $k = 4$ implies $N \approx (Z_2)^2$, $|X| = 4$.
(4) $k \neq 5$.

□

**Theorem 3.4([5], R. C. Bose):** There exists a finite projective plane of order $n$ if and only if there exists a complete set of $(n-1)$ MOLS (Mutually Orthogonal Latin Squares) of size $n$.

□

We give below a result about the nature of the number $n$ for the existence of complete set of MOLS of that size.

**Theorem 3.5** ([6], page 63)**:** If $n$ is prime or prime power then there exist $(n-1)$ MOLS (Mutually Orthogonal Latin Squares) of size $n$.

□

Actually, if there exists finite field formed by $n$ elements, namely, the elements {0, 1, 2, …., n-1} then there exist (n-1) MOLS (**but not conversely**) and it is known that such field exists precisely when $n$ is prime or power of prime. This paper is essentially about the important unsolved conjecture: There exists a full set of (*n*-1) MOLS if and only if *n* is prime or a power of prime.

4. **Existence of Finite Projective Planes:** The above mentioned conjecture will be done by proving converse of the above given Theorem 3.5, i.e. one needs to show that if (*n*-1) MOLS exist for some (given) *n* then that *n* must be either prime or prime power.

    Note that the existence of *fpp* of order *n* implies the existence of a subset of $S_n$ made up of $n(n-1)$ permutations which form a set of sharply 2-transitive permutations.
    
    **Now,**

    (A1) If these permutations already form a group of order $n(n-1)$ containing a regular normal subgroup of order *n* then *n* must be a prime power by stated equivalent theorems (Theorem 3.1, Theorem 3.2 and Theorem 3.3).

    (A2) If these permutations do not contain identity permutation then this set is not a group. But in this case, can be seen that it is actually a **coset.** From this coset the group with desired properties can be obtained by the following simple operation on this coset.

    (a) We pick a permutation in this set and operate on it by its inverse (in $S_n$) then the picked up permutation will change to the identity.

    (b) We then operate this same permutation (i.e. the inverse of the picked up permutation) on all the other permutations in this set of $n(n-1)$ permutations. We can see that we have a new set of sharply 2-transitive permutations containing identity permutation. This set is actually the group of sharply 2-transitive permutations, *G* say, containing a regularly normal subgroup, *N* say, formed by $(n-1)$ fixed-point-free permutations of degree *n* and identity.

    **Some Simple Observations:**

    (i) By elementary counting it can be seen that any set of sharply 2-transitive permutations contains just $n(n-1)$ sharply 2-transitive permutations by which the pairs of *n* elements get permuted

transitively and in which at most $(n-1)$ fixed-point-free permutations of degree $n$ can be present.

(ii) Now, since identity is present in the new set there must be present $(n-1)$ fixed-point-free permutations of degree $n$ in this set (to avoid columns with repeated entries) in order to form the Latin Square containing transformed row corresponding to identity element and the transformed rows corresponding to $(n-1)$ fixed-point-free elements, assured by the existence of *fpp*.

(iii) It is easy to check that $(n-1)$ fixed-point-free elements together with identity form a closed subset of this new set and so a subgroup.

(iv) The new set containing a given set of $(n-1)$ fixed-point-free elements is unique and contains exactly $n(n-1)$ elements, so the product of any two elements in this new set must belong to this new set and so, being closed, the new set also must be a group.

(v) It is easy to see that the $(n-1)$ fixed-point-free permutations are conjugate in $G$ and for any two distinct symbols among the $n$ symbols there is a unique fixed-point-free permutation which takes one symbol to the other.

(vi) These fixed-point-free elements being conjugate in $G$ will have same order $p$, where $p$ is a prime dividing $n$, and being transitive, will form together with identity a regularly normal subgroup $N$ of $G$.

Thus, when the set of 2-transitive $n(n-1)$ permutations, assured by the existence of a projective plane, do not contain identity permutation and so not a group then it is a **coset** from which we can obtain the sharply 2-transitive group of order $n(n-1)$ containing a regularly normal subgroup of order $n$, made up of $(n-1)$ fixed-point-free elements and the identity element, by the action of the inverse of some picked up permutation in this coset on the entire coset. And therefore again $n$ **must be a prime power** by stated equivalent theorems (Theorem 3.1, Theorem 3.2 and Theorem 3.3).

**Illustrative Example:** For the *fpp* of order $n = 4$ we saw above that the group of order $n(n-1)$ was $A_4$, the alternating group on four symbols, containing a normal subgroup of order $n = 4$ made up of three permutations without containing a fixed point and identity. But instead of taking the elements of this group if we take the other elements of $S_4$ then they do not form a group but they together also satisfy the property described in section 4 and can be used to form 3 MOLS!

| Permutation(disjoint cycle form) | Transformed row |
|---|---|
| (1234) | 2 3 4 1 |
| (1432) | 4 1 2 3 |
| (1324) | 3 4 2 1 |
| (1342) | 3 1 4 2 |
| (1243) | 2 4 1 3 |
| (1423) | 4 3 1 2 |
| (1)(2)(34) | 1 2 4 3 |
| (1)(3)(24) | 1 4 3 2 |
| (1)(4)(23) | 1 3 2 4 |
| (2)(3)(14) | 4 2 3 1 |
| (2)(4)(13) | 3 2 1 4 |
| (3)(4)(12) | 2 1 3 4 |

The 3 MOLS are:

```
1 3 2 4        1 4 3 2        1 2 4 3
2 4 1 3        2 3 4 1        2 1 3 4
3 1 4 2        3 2 1 4        3 4 2 1
4 2 3 1        4 1 2 3        4 3 1 2
```

Now, apply the inverse permutation of the first permutation, namely, the cycle (1234) [which is (1432), since (1432)(1234) = (1)(2)(3)(4)] on every permutation represented in terms of cycle decomposition in the first column of above 2-columned table. One can see that the permutations corresponding to the group $A_4$ (as desired) are produced!!

We now proceed with

**Alternative Way of Implying existence of Group from existence of ($n$-1) MOLS:**

How to obtain the existence of a sharply 2-transitive subgroup of $S_n$ of order $n(n-1)$ containing a regularly normal subgroup of order

$n$, from the existence of ($n$-1) MOLS? We now discuss it in brief and demonstrate this for the above considered case: $n = 4$.

Latin Squares may be thought of as matrices of size $n$ whose first row can be taken as [0 1 2 3 ……… $n$-1] and the other rows as suitable permutations of the elements of this row so that all rows will automatically contain each element of this set only once and by properly selecting these permutations the care can be taken of every column to contain each element of this set only once as desired by the formal definition of the Latin Square. Now suppose we have ($n$-1) MOLS made up of elements from the set $S$ having first row as [0 1 2 3 ……… $n$-1], i.e. the first row of all MOLS is the row vector [0 1 2 3 ……… $n$-1]. Now, it is easy to see that if we will apply any identical permutation that permutes rows of a Latin Square on all the given ($n$-1) MOLS simultaneously then the new set of matrices that results will also forms a set of ($n$-1) MOLS.

So we carry out the following **Steps** to convert given set of MOLS, say $\{L_1^*, L_2^*, L_3^*, \cdots, L_{n-1}^*\}$ to another (equivalent) new set of MOLS, say $\{L_1, L_2, L_3, \cdots, L_{n-1}\}$ so that the permutations in the form of transformed rows in new set $\{L_i^C + 1, 1 \leq i \leq n-1\}$ will form a sharply 2-transitive group of order $n(n-1)$ containing a normal subgroup of order $n$ formed by $(n-1)$ fixed-point-free permutations and identity.

(1) Let $\{L_1^*, L_2^*, L_3^*, \cdots, L_{n-1}^*\}$ be the set of ($n$-1) MOLS that exist for given order $n$.

$$L_1^* = \begin{bmatrix} 0 & 1 & 2 & . & . & n-1 \\ . & . & . & . & . & . \\ . & . & . & . & . & . \\ . & . & . & . & . & . \\ . & . & . & . & . & . \\ . & . & . & . & . & . \end{bmatrix}, L_2^* = \begin{bmatrix} 0 & 1 & 2 & . & . & n-1 \\ . & . & . & . & . & . \\ . & . & . & . & . & . \\ . & . & . & . & . & . \\ . & . & . & . & . & . \\ . & . & . & . & . & . \end{bmatrix}, \cdots,$$

(2) Apply a suitably chosen same permutation, permuting the rows of Latin Squares, identically and simultaneously on each of the Latin

Squares in given set of (*n*-1) MOLS, so that the new set of (*n*-1) MOLS that results will contain a Latin Square, say $L_1$ having its first row as $[0 \ 1 \ 2 \ . \ . \ n-1]$ and also its first column as

$$\begin{bmatrix} 0 \\ 1 \\ 2 \\ . \\ . \\ n-1 \end{bmatrix}$$

Thus, the **new set** of MOLS will be as follows:

$$L_1 = \begin{bmatrix} 0 & 1 & 2 & . & . & n-1 \\ 1 & . & . & . & . & . \\ 2 & . & . & . & . & . \\ . & . & . & . & . & . \\ n-1 & . & . & . & . & . \end{bmatrix}, \cdots, L_{n-1} = \begin{bmatrix} 0 & 1 & 2 & . & . & n-1 \\ . & . & . & . & . & . \\ . & . & . & . & . & . \\ . & . & . & . & . & . \\ . & . & . & . & . & . \end{bmatrix}$$

(3) Let $\{L_1, L_2, L_3, \cdots, L_{n-1}\}$ the above obtained **new set** of MOLS. Find transpose of each of the Latin Square matrix in the above newly obtained set of Latin Square matrices and add unit in each entry of these matrices. In other words, in each transpose matrix add *n* by *n* matrix with each of its entry equal to 1. Thus, we get the new set of matrices as follows:

$$L_1^C + 1 = \begin{bmatrix} 1 & 2 & 3 & . & . & n \\ 2 & . & . & . & . & . \\ 3 & . & . & . & . & . \\ . & . & . & . & . & . \\ n & . & . & . & . & . \end{bmatrix}, \cdots, L_{n-1}^C + 1 = \begin{bmatrix} 1 & . & . & . & . & . \\ 2 & . & . & . & . & . \\ 3 & . & . & . & . & . \\ . & . & . & . & . & . \\ n & . & . & . & . & . \end{bmatrix}$$

(4) Each matrix in this new set of matrices $\{L_i^C + 1, 1 \leq i \leq n-1\}$ is also a Latin Square but these matrices have additional structure, in the sense that the rows of $\{L_i^C + 1, 1 \leq i \leq n-1\}$ correspond to transformed rows of permutations. Since each matrix has $n$ rows in it and in all there are ($n$-1) such matrices therefore there are in all $n(n-1)$ transformed rows of permutations. It can be verified that for those $n$ for which MOLS $\{L_i, 1 \leq i \leq n-1\}$ exist the permutations represented by the transformed rows in the corresponding $\{L_i^C + 1, 1 \leq i \leq n-1\}$ are sharply 2-transitive and they together form the group of sharply 2-transitive permutations of order $n(n-1)$. Further it contains a normal subgroup of order $n$ formed by fixed point free permutations and identity. In the matrix $L_1^C + 1$ the first (transformed) row (of permutation) corresponds to identity permutation with cycle structure (1)(2)(3)....(n). The other rows of this same matrix represent ($n$-1) fixed point free permutations for if any permutation among these ($n$-1) permutations is not fixed point free then some column of $L_1^C + 1$ will contain some entry twice which is forbidden by definition.

**Example:** Let us consider the case for $n = 4$: Let $\{L_1^{*C} + 1, L_2^{*C} + 1, L_3^{*C} + 1\}$ be the given MOLS formed by the transformed rows of permutations which are

```
1 3 2 4        1 4 3 2        1 2 4 3
2 4 1 3        2 3 4 1        2 1 3 4
3 1 4 2        3 2 1 4        3 4 2 1
4 2 3 1        4 1 2 3        4 3 1 2
```

These permutations in the form of transformed rows when expressed as product of disjoint cycles happen to be

```
(1)(23)(4)     (1)(24)(3)     (1)(2)(34)
(1243)         (1234)         (12)(3)(4)
(1342)         (13)(2)(4)     (1324)
(14)(2)(3)     (1432)         (1423)
```

These permutations together form a set of sharply 2-transitive permutations but they together do not form a group.

Now, by taking transpose and subtracting from it the matrix of same size having all its elements equal to unity from the Latin squares given above formed by the transformed rows we get three MOLS in the standard form, $\{L_1^*, L_2^*, L_3^*\}$, as follows:

$$
\begin{array}{ccc}
0\ 1\ 2\ 3 & 0\ 1\ 2\ 3 & 0\ 1\ 2\ 3 \\
2\ 3\ 0\ 1 & 3\ 2\ 1\ 0 & 1\ 0\ 3\ 2 \\
1\ 0\ 3\ 2 & 2\ 3\ 0\ 1 & 3\ 2\ 1\ 0 \\
3\ 2\ 1\ 0 & 1\ 0\ 3\ 2 & 2\ 3\ 0\ 1 \\
\end{array}
$$

Now, we apply suitable permutation (that permutes the rows) identically and simultaneously on all $\{L_1^*, L_2^*, L_3^*\}$ given above, <u>in the present case we choose permutation that interchanges second row with third row everywhere</u> and obtain $\{L_1, L_2, L_3\}$ as follows:

$$
\begin{array}{ccc}
0\ 1\ 2\ 3 & 0\ 1\ 2\ 3 & 0\ 1\ 2\ 3 \\
1\ 0\ 3\ 2 & 2\ 3\ 0\ 1 & 3\ 2\ 1\ 0 \\
2\ 3\ 0\ 1 & 3\ 2\ 1\ 0 & 1\ 0\ 3\ 2 \\
3\ 2\ 1\ 0 & 1\ 0\ 3\ 2 & 2\ 3\ 0\ 1 \\
\end{array}
$$

From these $\{L_1, L_2, L_3\}$ we now obtain $\{L_1^C + 1, L_2^C + 1, L_3^C + 1\}$ as follows:

$$
\begin{array}{ccc}
1\ 2\ 3\ 4 & 1\ 3\ 4\ 2 & 1\ 4\ 2\ 3 \\
2\ 1\ 4\ 3 & 2\ 4\ 3\ 1 & 2\ 3\ 1\ 4 \\
3\ 4\ 1\ 2 & 3\ 1\ 2\ 4 & 3\ 2\ 4\ 1 \\
4\ 3\ 2\ 1 & 4\ 2\ 1\ 3 & 4\ 1\ 3\ 2 \\
\end{array}
$$

These permutations in the form of transformed rows when expressed as product of disjoint cycles look like

$$
\begin{array}{ccc}
(1)(2)(3)(4) & (1)(234) & (1)(243) \\
(12)(34) & (124)(3) & (123)(4) \\
(13)(24) & (132)(4) & (134)(2) \\
(14)(23) & (143)(2) & (142)(3) \\
\end{array}
$$

Clearly, these permutations together form a group as expected of sharply 2-transitive permutations and also they together form alternating group, $A_4$. Further, the permutations corresponding to transformed rows in $L_1^C + 1$ form regularly normal subgroup of $A_4$ of order four formed by three fixed point free permutations and identity. The above can also be explained as follows: $\{L_1^C + 1, L_2^C + 1, L_3^C + 1\}$ can be obtained from $\{L_1^{*C} + 1, L_2^{*C} + 1, L_3^{*C} + 1\}$ by interchanging columns 2 and 3 in $\{L_1^{*C} + 1, L_2^{*C} + 1, L_3^{*C} + 1\}$ and this was aimed at arranging first row of $L_1^{*C} + 1$ in increasing order.

5. **More on MOLS:** In this section we put forward some more observations about MOLS. Formally, a **Latin Square** (of order $n$), $L = (a_{ij})$, is **defined** as a square matrix containing in all $n^2$ elements, and these elements are taken from a set containing $n$ elements, whose each row and column contains every element from the chosen set of $n$ elements exactly once. Two Latin Squares, $L = (a_{ij})$ and $M = (b_{ij})$ are said to be **mutually orthogonal** if the $n^2$ ordered pairs $(a_{ij}, b_{ij})$ are all different. Let us take the set of $n$ elements, $S$, as $S = \{0,1,2,3,\cdots,n-1\}$. All pairs $(a_{ij}, b_{ij})$ are distinct meant that every possible pair formed by taking product $S \times S$ has appeared and appeared only once. Suppose $\{L_1, L_2, L_3, \cdots, L_{n-1}\}$ be a set of Latin Squares. Let us denote $L_k = (a_{ij}^k)$. Let us form ordered $(n-1)$-tuples, formed by taking direct product $L_1 \times L_2 \times L_3 \times \cdots \times L_{n-1}$. It will constitute a set of $(n-1)$-tuples $(a_{ij}^1, a_{ij}^2, a_{ij}^3, \cdots, a_{ij}^{n-1})$. There will be in all $n^2$ such ordered $(n-1)$-tuples, as by taking $(ij)$-th entry from each Latin Square we will form one $(n-1)$-tuple, and there are in all $n^2$ such entries. The above definition about mutual orthogonality of two Latin Squares can be easily extended as follows: The Latin Squares $\{L_1, L_2, L_3, \cdots, L_{n-1}\}$ will be mutually orthogonal, i.e. they together are $(n-1)$ MOLS if and only if these $(n-1)$-tuples $(a_{ij}^1, a_{ij}^2, a_{ij}^3, \cdots, a_{ij}^{n-1})$ are all different. These $(n-1)$-tuples are different meant that if we consider any two such $(n-1)$-tuples then either all their coordinates are different or they intersect at most at one coordinate, i.e. let $(a_{ij}^1, a_{ij}^2, a_{ij}^3, \cdots, a_{ij}^{n-1})$ and $(b_{ij}^1, b_{ij}^2, b_{ij}^3, \cdots, b_{ij}^{n-1})$ be any two $(n-1)$-tuples then either $a_{ij}^r \neq b_{ij}^r$ for all $r$, or $a_{ij}^r = b_{ij}^r$ for at most some one $r$. Check that these tuples taken together are such that every $a_{ij}^r, 1 \leq r \leq n-1$ takes every value from set $S = \{0,1,2,3,\cdots,n-1\}$ exactly $n$ times.

**Proposition 5.1:** Every Latin Square can be thought of as matrix of size $n$ which contains $n$ distinct determinantal monomials or which is made up of $n$ distinct determinantal monomials:

$$a_{1\sigma^1(1)} a_{2\sigma^1(2)} a_{3\sigma^1(3)} \cdots a_{n\sigma^1(n)},$$

$$a_{1\sigma^2(1)} a_{2\sigma^2(2)} a_{3\sigma^2(3)} \cdots a_{n\sigma^2(n)},$$

$$\vdots$$

$$a_{1\sigma^n(1)} a_{2\sigma^n(2)} a_{3\sigma^n(3)} \cdots a_{n\sigma^n(n)}.$$

where

$\sigma^k$ are $n$ distinct permutations such that $a_{i\sigma^k(i)} = k-1$, and $1 \leq i \leq n$ and $1 \leq k \leq n$.

**Proof:** First monomial ensures that number zero ("0") appears in every row (and column) in a matrix only once. Second monomial ensures that number one ("1") appears in every row (and column) in the same matrix only once. Continuing on these lines the last ($n$-th) monomial ensures that the last number ("$n$-1") appears in every row (and column) in the same matrix only once. In effect, every number from set $S = \{0,1,2,3,\cdots,n-1\}$ appears in every row/column of the same matrix exactly one times. Hence, the matrix must be a Latin Square.

This proposition implies that a Latin Square of order *n* can be thought of as matrix which contains matrix elements from set $S = \{0,1,2,3,\cdots,n-1\}$ and these elements are placed in such a way that there are *n* distinct monomials, one made up of all elements equal to '0', other of '1', …., other of '(n-1)'.

**Proposition 5.2:** Two Latin Squares are mutually orthogonal if and only if any two determinantal monomials, one chosen from *n* monomials corresponding to first Latin Square and the other chosen from *n* monomials corresponding to second Latin Square, have at most one intersection.

**Proof:** By definition two Latin Squares are mutually orthogonal if the $n^2$ ordered pairs $(a_{ij}, b_{ij})$ are all different. Let

$$a_{1\sigma^i(1)} a_{2\sigma^i(2)} a_{3\sigma^i(3)} \cdots a_{n\sigma^i(n)}, \quad b_{1\sigma^j(1)} b_{2\sigma^j(2)} b_{3\sigma^j(3)} \cdots b_{n\sigma^j(n)}$$

be the two determinantal monomials and let $a_{m\sigma^i(m)} = p, b_{m\sigma^j(m)} = q$ and $a_{n\sigma^i(n)} = p, b_{n\sigma^j(n)} = q$, i.e. let these monomials intersect at two places, then when two Latin Squares will be superimposed and pairs will be formed we will have two identical pairs, $(p,q)$, which contradicts orthogonality. The converse is also straightforward.
□

**Example:** Let *n* = 4. There will be in all three MOLS, $L_1, L_2, L_3$. The determinantal monomials corresponding to $L_1$ are

$$a_{11}a_{22}a_{33}a_{44},$$
$$a_{12}a_{21}a_{34}a_{43},$$
$$a_{13}a_{24}a_{31}a_{42},$$
$$a_{14}a_{23}a_{32}a_{41}$$

The determinantal monomials corresponding to $L_2$ are

$$a_{11}a_{23}a_{34}a_{42},$$

$$a_{12}a_{24}a_{33}a_{41},$$

$$a_{13}a_{21}a_{32}a_{44},$$

$$a_{14}a_{22}a_{31}a_{43}$$

The determinantal monomials corresponding to $L_3$ are

$$a_{11}a_{24}a_{32}a_{43},$$

$$a_{12}a_{23}a_{31}a_{44},$$

$$a_{13}a_{22}a_{34}a_{41},$$

$$a_{14}a_{21}a_{33}a_{42}$$

At the positions depicted by **first** monomials in the determinantal monomials corresponding to $L_1, L_2, L_3$ we substitute **'0'**. Similarly, at the positions depicted by **second, third and fourth** monomials in the determinantal monomials corresponding to $L_1, L_2, L_3$ we substitute **'1','2', and '3'** respectively. This leads us to MOLS $L_1, L_2, L_3$ as shown below:

$$L_1 = \begin{bmatrix} 0 & 1 & 2 & 3 \\ 1 & 0 & 3 & 2 \\ 2 & 3 & 0 & 1 \\ 3 & 2 & 1 & 0 \end{bmatrix},$$

$$L_2 = \begin{bmatrix} 0 & 1 & 2 & 3 \\ 2 & 3 & 0 & 1 \\ 3 & 2 & 1 & 0 \\ 1 & 0 & 3 & 2 \end{bmatrix},$$

$$L_3 = \begin{bmatrix} 0 & 1 & 2 & 3 \\ 3 & 2 & 1 & 0 \\ 1 & 0 & 3 & 2 \\ 2 & 3 & 0 & 1 \end{bmatrix}$$

**Observation 5.1:** Let $\{L_1, L_2, L_3, \cdots, L_{n-1}\}$ be a set of (n-1) MOLS and

$$L_1 = \begin{bmatrix} 0 & 1 & 2 & . & . & n-1 \\ 1 & . & . & . & . & . \\ 2 & . & . & . & . & . \\ . & . & . & . & . & . \\ n-1 & . & . & . & . & . \end{bmatrix}$$

then there exists a cyclic group, $G$, such that the other MOLS $\{L_2, L_3, \cdots, L_{n-1}\}$ can be obtained from $L_1$ by action of this cyclic group. The elements of this group permute rows of $L_1$ except the first row. Let $g$ be the generator of this group, then $g(L_1) = L_2$, $g^2(L_1) = L_3$, $g^3(L_1) = L_4$, ...., $g^{n-1}(L_1) = L_1$, i.e. $g^{n-1} = e$, the identity. Further, the set of MOLS $\{L_1, L_2, L_3, \cdots, L_{n-1}\}$ is closed under the action of this group, i.e. if we operate on any Latin Square $L_i$ some element, $\alpha \in G$ then we always get some other Latin Square $L_j$, i.e. $\alpha(L_i) = L_j$.

□

**Example 1:** Consider MOLS $L_1, L_2, L_3$ given on last page. For the sake of convenience we call first row of Latin Squares as zero[th] row, second row as first row, and so on. We have $G = \{e, \alpha, \alpha^2\}$, where $\alpha = (0)(123)$. Clearly, $\alpha(L_1) = L_2$, $\alpha^2(L_1) = \alpha(L_2) = L_3$.

□

**Example 2:** Let us take $n = 9$. $\{L_1, L_2, L_3, \cdots, L_8\}$ be a set of 8 MOLS and

$$L_1 = \begin{bmatrix} 0 & 1 & 2 & 3 & 4 & 5 & 6 & 7 & 8 \\ 1 & 2 & 0 & 4 & 5 & 3 & 7 & 8 & 6 \\ 2 & 0 & 1 & 5 & 3 & 4 & 8 & 6 & 7 \\ 3 & 4 & 5 & 6 & 7 & 8 & 0 & 1 & 2 \\ 4 & 5 & 3 & 7 & 8 & 6 & 1 & 2 & 0 \\ 5 & 3 & 4 & 8 & 6 & 7 & 2 & 0 & 1 \\ 6 & 7 & 8 & 0 & 1 & 2 & 3 & 4 & 5 \\ 7 & 8 & 6 & 1 & 2 & 0 & 4 & 5 & 3 \\ 8 & 6 & 7 & 2 & 0 & 1 & 5 & 3 & 4 \end{bmatrix}$$

For the sake of convenience we call first row of Latin Squares as zero$^{th}$ row, second row as first row, and so on. We have $G = \{e, \alpha, \alpha^2, \cdots, \alpha^7\}$, where $\alpha = (0)(13782654)$, $\alpha^2 = (0)(1725)(3864)$, $\alpha^3 = (0)(18532476)$, $\alpha^4 = (0)(12)(36)(48)(57)$, $\alpha^5 = (0)(16742358)$, $\alpha^6 = (0)(1527)(3468)$, $\alpha^7 = (0)(14562873)$, $\alpha^8 = (0)(1)(2)(3)(4)(5)(6)(7)(8) = e$.

It is easy to check that we get $\alpha(L_1) = L_2$, where

$$L_2 = \begin{bmatrix} 0 & 1 & 2 & 3 & 4 & 5 & 6 & 7 & 8 \\ 3 & 4 & 5 & 6 & 7 & 8 & 0 & 1 & 2 \\ 6 & 7 & 8 & 0 & 1 & 2 & 3 & 4 & 5 \\ 7 & 8 & 6 & 1 & 2 & 0 & 4 & 5 & 3 \\ 1 & 2 & 0 & 4 & 5 & 3 & 7 & 8 & 6 \\ 4 & 5 & 3 & 7 & 8 & 6 & 1 & 2 & 0 \\ 5 & 3 & 4 & 8 & 6 & 7 & 2 & 0 & 1 \\ 8 & 6 & 7 & 2 & 0 & 1 & 5 & 3 & 4 \\ 2 & 0 & 1 & 5 & 3 & 4 & 8 & 6 & 7 \end{bmatrix}$$

We can easily check that we get $\alpha^2(L_1) = \alpha(L_2) = L_3$,
$\alpha^3(L_1) = \alpha^2(L_2) = \alpha(L_3) = L_4$, $\alpha^4(L_1) = \alpha^3(L_2) = \alpha^2(L_3) = \alpha(L_4) = L_5$,
$\alpha^5(L_1) = \alpha^4(L_2) = \alpha^3(L_3) = \alpha^2(L_4) = \alpha(L_5) = L_6$,
$\alpha^6(L_1) = \alpha^5(L_2) = \alpha^4(L_3) = \alpha^3(L_4) = \alpha^2(L_5) = \alpha(L_6) = L_7$,
$\alpha^7(L_1) = \alpha^6(L_2) = \alpha^5(L_3) = \alpha^4(L_4) = \alpha^3(L_5) = \alpha^2(L_6) = \alpha(L_7) = L_8$ and
finally, $\alpha^8(L_1) = L_1$.

□

**Some Important Remarks:** When order $n$ is prime, $p$ say, then there is no divisor for this order other than 1 or $n = p$. So, if we consider sequence $i_1 i_2 i_3 \cdots i_p$ where $i_1 = r, 1 \le r \le p$ and $|i_{r+1} - i_r| = k, 1 \le k \le p-1$ then this sequence will contain all numbers from 1 to $p$ and all these sequences will be different in the sense that any two numbers in them will obey 2-transitivity, i.e. there will be in all $n(n-1)$ such sequences (forming transformed rows of permutations) and when all these sequences are kept one below the other then in any two columns we will be having every pair among the pairs (1, 2), (1, 3), **....**, (1, n), (2, 1), (2, 3), **....**, (2, n), (3, 1), (3,2), (3,4), **....**, ((n-1), n), (n ,1), (n, 2), **....**, (n, (n-1)) and it will be appearing only once. Using this fact we can construct $\{L_1^C + 1, L_2^C + 1, L_3^C + 1, \cdots, L_{n-1}^C + 1\}$ as follows. We also can see that the rows of $L_1^C + 1$ are the transformed rows of permutations that contains identity and $(p-1)$ fixed point free permutations which are all $p$-cycles. Further these permutations form the subgroup of order $p$ in the group of order $(p)(p-1)$ that is formed when permutations corresponding to transformed rows of all $\{L_1^C + 1, L_2^C + 1, L_3^C + 1, \cdots, L_{n-1}^C + 1\}$ are taken together. Following are

$$L_1^C + 1 = \begin{bmatrix} 1 & 2 & 3 & \cdots & p \\ 2 & 3 & 4 & \cdots & 1 \\ 3 & 4 & 5 & \cdots & 2 \\ \vdots & \vdots & \vdots & \cdots & \\ p & 1 & 2 & \cdots & p-1 \end{bmatrix}$$

$$L_2^C + 1 = \begin{bmatrix} 1 & 3 & 5 & \cdots & p-1 \\ 2 & 4 & 6 & \cdots & p \\ 3 & 5 & 7 & \cdots & 1 \\ \vdots & \vdots & \vdots & \cdots & \vdots \\ p & 2 & 4 & \cdots & p-2 \end{bmatrix}$$

$$L_3^C + 1 = \begin{bmatrix} 1 & 4 & 7 & \cdots & p-2 \\ 2 & 5 & 8 & \cdots & p-1 \\ 3 & 6 & 9 & \cdots & p \\ \vdots & \vdots & \vdots & \cdots & \vdots \\ p & 3 & 6 & \cdots & p-3 \end{bmatrix}$$

$$\vdots$$

$$L_{n-1}^C + 1 = L_{p-1}^C + 1 = \begin{bmatrix} 1 & p & p-1 & \cdots & 2 \\ 2 & 1 & p & \cdots & 3 \\ 3 & 2 & 1 & \cdots & 4 \\ \vdots & \vdots & \vdots & \cdots & \vdots \\ p & p-1 & p-2 & \cdots & 1 \end{bmatrix}$$

All the rows together in $\{L_1^C + 1, L_2^C + 1, L_3^C + 1, \cdots, L_{n-1}^C + 1\}$ are transformed rows for permutations and form the set of in all $n(n-1)$ sharply 2-transitive permutations and actually they form a group. Further, all rows together of $L_1^C + 1$ form a subgroup (of the group of order $n(n-1)$ just mentioned) and is made up of identity and $(n-1)$ fixed point free permutations which are $n$-cycles. The first row of $L_1^C + 1$ corresponds to transformed row for permutation identity, $e = (1)(2)(3)\cdots(p)$. The second row corresponds to transformed row for permutation which is $p$-cycle, $(123\cdots p)$. The third row corresponds to transformed row for permutation which is *p*-cycle, $(1357\cdots(p-1))$. The fourth row corresponds to transformed row for permutation which is *p*-cycle, $(147\cdots(p-2))$. The last row corresponds to transformed row for permutation which is *p*-cycle, $(1p(p-1)(p-2)\cdots 2)$. It is easy to check that if we take into consideration the entries in some row in the adjacent (neighboring) columns of $L_j^C + 1$ then they are $k$ and $(k+j)\mod(p)$ respectively. If we arrange all the rows in the matrices $\{L_1^C + 1, L_2^C + 1, L_3^C + 1, \cdots, L_{n-1}^C + 1\}$ one below the other and form a bigger matrix of size *n* by $n(n-1)$ and consider any two columns of this matrix then there is exactly one time appearance of every pair among the pairs (1, 2), (1, 3), ...., (1, n), (2, 1), (2, 3), ...., (2, n), (3, 1), (3,2), (3,4), ...., ((n-1), n), (n ,1), (n, 2), ...., (n, (n-1)) in each such pair of columns $(i, j)$ such that $1 \le i \le n$ and $1 \le j \le n$.

When order *n* is prime power, $p^k, k>1$, say, then still we can form the matrices $\{L_1^C + 1, L_2^C + 1, L_3^C + 1, \cdots, L_{n-1}^C + 1\}$ which possess the same properties as above when *n* was only prime and was not a prime power. All the rows together in $\{L_1^C + 1, L_2^C + 1, L_3^C + 1, \cdots, L_{n-1}^C + 1\}$ form the set of $n(n-1)$ sharply 2-transitive permutations and actually they form a group. Further, we will see that in this case also all rows together of $L_1^C + 1$ form a subgroup (of the group of order $n(n-1)$ just mentioned) and is made up of identity and $(n-1)$ fixed point free permutations which are <u>product of $p^{k-1}$ number of totally disjoint *p*-cycles.</u> The first row of $L_1^C + 1$ corresponds to transformed row for permutation identity, $e = (1)(2)(3)\cdots(p^k)$. The second row corresponds to transformed row for permutation which is product of $p^{k-1}$ number of totally disjoint *p*-cycles, $(123\cdots p)((p+1)\cdots 2p)((2p+1)\cdots 3p)\cdots((p^{k-1}-1)p+1)\cdots p^k)$. (Disjoint *p*-cycles as given above in the product are called <u>totally disjoint *p*-cycles</u> because they have no entry in common). The third and other rows corresponds to transformed row for permutation which are all of the same type, i.e. product of $p^{k-1}$ number of totally disjoint *p*-cycles and they are formed in such a way that these permutations represented by product of totally disjoint cycles and identity together <u>form a group of order $n$ which ultimately form subgroup of the group of order $n(n-1)$ made up of sharply 2-transitive permutations</u>. There are in all $(n-1)$ matrices, $L_j^C +1, 1 \leq j \leq (n-1)$, and $n$ rows of each matrix gives rise to $n$ sharply 2-transitive permutations thus together there are required $n(n-1)$ sharply 2-transitive permutations and they together form a group.

**Example 1:** Let $n = 5$. We construct $\{L_1^C +1, L_2^C +1, L_3^C +1, L_4^C +1\}$ as seen above for prime case:

$$L_1^C + 1 = \begin{bmatrix} 1 & 2 & 3 & 4 & 5 \\ 2 & 3 & 4 & 5 & 1 \\ 3 & 4 & 5 & 1 & 2 \\ 4 & 5 & 1 & 2 & 3 \\ 5 & 1 & 2 & 3 & 4 \end{bmatrix}$$

$$L_2^C + 1 = \begin{bmatrix} 1 & 3 & 5 & 2 & 4 \\ 2 & 4 & 1 & 3 & 5 \\ 3 & 5 & 2 & 4 & 1 \\ 4 & 1 & 3 & 5 & 2 \\ 5 & 2 & 4 & 1 & 3 \end{bmatrix}$$

$$L_3^C + 1 = \begin{bmatrix} 1 & 4 & 2 & 5 & 3 \\ 2 & 5 & 3 & 1 & 4 \\ 3 & 1 & 4 & 2 & 5 \\ 4 & 2 & 5 & 3 & 1 \\ 5 & 3 & 1 & 4 & 2 \end{bmatrix}$$

$$L_4^C + 1 = \begin{bmatrix} 1 & 5 & 4 & 3 & 2 \\ 2 & 1 & 5 & 4 & 3 \\ 3 & 2 & 1 & 5 & 4 \\ 4 & 3 & 2 & 1 & 5 \\ 5 & 4 & 3 & 2 & 1 \end{bmatrix}$$

It is easy to see that rows of $L_1^C + 1$ represent following permutations: $\{(1)(2)(3)(4)(5),(12345),(13524),(14253),(15432)\}$ and they together form group which is normal subgroup of group of permutations given in the second column of the two columned table given below.

| Matrix | Permutations |
|---|---|
| $L_1^C + 1$ | $\{(1)(2)(3)(4)(5),(12345),(13524),(14253),(15432)\}$ |
| $L_2^C + 1$ | $\{(1)(2354),(5)(1243),(4)(1325),(3)(1452),(2)(1534)\}$ |
| $L_3^C + 1$ | $\{(1)(2453),(3)(1254),(5)(1342),(2)(1435),(4)(1523)\}$ |
| $L_4^C + 1$ | $\{1)(25)(34),(4)(12)(35),(2)(13)(45),(5)(14)(23),(3)(15)(24)$ |

**Example 2:** We give below permutations obtained using rows of $L_1^C + 1$ (which represent transformed rows of permutations) forming group order $n$ which is normal subgroup of bigger group of order $n(n-1)$ formed by sharply 2-transitive permutations for $n = 2^2, 2^3, 3^2$ in the two columned table below:

| $n$ | Permutation |
|---|---|
| $2^2 = 4$ | {Identity,(12)(34),(13)(24),(14)(23)} |
| $2^3 = 8$ | {Identity, (12)(34)(56)(78), (13)(24)(57)(68), (14)(23)(58)(67), (15)(26)(37)(48), ((16)(25)(38)(47), (17)(28)(35)(46), (18)(27)(36)(45)} |
| $3^2 = 9$ | {identity, (123)(456)(789), (132)(465)(798), (147)(258)(369), (159)(267)(348), (168)(249)(357), (174)(285)(396), (186)(294)(375), (195)(276)(384)} |

# Acknowledgements

Thanks are due to Prof. Dr. M. R. Modak and Prof. Dr. S. A. Katre for some useful discussions.